\documentclass[12pt]{article}
\usepackage{latexsym,amsmath,amssymb,amsfonts,epsfig,graphicx,cite,psfrag}
\usepackage{eepic,color,colordvi,amscd,amsthm}

\newtheorem{theo}{Theorem}[section]

\newtheorem{coro}[theo]{Corollary}

\theoremstyle{definition}

\def\qed{\hfill \rule{4pt}{7pt}}
\def\pf{\noindent {\it Proof.} }

\makeatletter
\renewcommand \theequation  {\ifnum \c@chapter>\z@ \thechapter.\fi
\ifnum \c@section>\z@ \@arabic\c@section.\fi \ifnum
\c@subsection>\z@ \@arabic\c@subsection.\fi\ifnum
\c@subsubsection>\z@
\@arabic\c@subsubsection.\fi\@arabic\c@equation}
\makeatother

\makeatletter
\@addtoreset{equation}{section}
\@addtoreset{figure}{section}
\makeatother

\renewcommand{\theequation}{\arabic{section}.\arabic{equation}}

\begin{document}
\begin{center}
{\bf\large Linked Partitions and Permutation Tableaux}\\ \vskip 6mm

{\small William Y.C. Chen$^1$, Lewis H. Liu$^2$, Carol J. Wang$^3$\\[2mm]

$^1$Center for Applied Mathematics\\
Tianjin University\\ Tianjin 300072, P.R. China\\[3mm]

$^2$Center for Combinatorics, LPMC-TJKLC\\
Nankai University\\
 Tianjin 300071, P.R. China \\[3mm]

$^3$Department of Mathematics \\
Beijing Technology and
Business University\\ Beijing 100048, P.R. China \\[3mm]
$^1$chenyc@tju.edu.cn, $^2$lewis@cfc.nankai.edu.cn, $^3$wang$\_$jian@th.btbu.edu.cn}

\end{center}

\begin{abstract}
Linked partitions are introduced by Dykema
 in the study of transforms in free probability
theory, whereas permutation tableaux are introduced by Steingr\'{i}msson and Williams in the study of totally positive
Grassmannian cells. Let $[n]=\{1,2,\ldots,n\}$.
Let $L(n,k)$ denote the set of  linked partitions of
  $[n]$ with $k$ blocks,  let $P(n,k)$ denote the set of  permutations of $[n]$ with $k$ descents, and  let $T(n,k)$ denote the set of permutation tableaux of length $n$ with $k$ rows. Steingr\'{i}msson and Williams found a bijection between the set of permutation tableaux of length $n$ with $k$ rows and the set of permutations of $[n]$ with $k$ weak excedances. Corteel and Nadeau gave a bijection from the set of permutation tableaux of length $n$ with $k$ columns to the set of permutations of $[n]$ with $k$ descents.
 In this paper, we establish a bijection
 between $L(n,k)$ and $P(n,k-1)$ and a bijection between $L(n,k)$ and $T(n,k)$. Restricting
 the latter bijection to  noncrossing linked partitions,
 we find that the corresponding  permutation tableaux
 can be characterized by pattern avoidance.

\end{abstract}



\vskip6mm
\section{Introduction}

The notion of linked partitions was introduced by Dykema \cite{Dyke07} in the study of the unsymmetrized T-transform in free probability theory.  Let $[n]=\{1,2,\ldots,n\}$. A linked partition of $[n]$ is a collection of nonempty subsets $B_1, B_2, \ldots, B_k$ of $[n]$, called blocks, such that the union of $B_1, B_2, \ldots, B_k$ is $[n]$ and any two distinct blocks are nearly disjoint. Two blocks $B_i$ and $B_j$  are said to be nearly disjoint if for any $k \in B_i\cap B_j$, one of the following conditions holds:
\begin{itemize}
\item[(a)] $k=\min(B_i), |B_i|>1$ and $k\neq \min(B_j)$, or
\item[(b)] $k=\min(B_j), |B_j|>1$ and $k\neq \min(B_i)$.
\end{itemize}

\begin{figure}[h]
\begin{center}
\setlength{\unitlength}{0.2cm}
\begin{picture}(40,9)
\multiput(0,3)(5, 0){9}{\circle{0.5}}\put(-.6,0){$1$}\put(4.6,0){$2$}
\put(9.6,0){$3$}\put(14.6,0){$4$}\put(19.6,0){$5$}\put(24.6,0){$6$}\put(29.6,0){$7$}
\put(34.6,0){$8$}\put(39.6,0){$9$}\qbezier(10,3)(25,11)(40,3)\qbezier(0,3)(7.5,11)(15,3)
\qbezier(5,3)(7.5,6)(10,3)\qbezier(20,3)(22.5, 6)(25,3)
\qbezier(25,3)(27.5,6)(30,3)\qbezier(0,3)(2.5,6)(5,3)
\end{picture}
\caption{The linear representation of a linked partition.}\label{LP}
\end{center}
\end{figure}
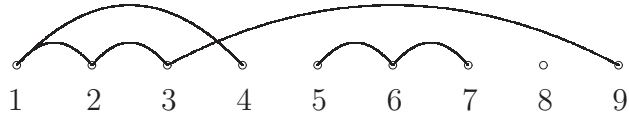

 The linear representation of a linked partition
  was introduced by Chen, Wu and Yan \cite{CWY08}.
For a linked partition $\tau$ of $[n]$, list the $n$ vertices in a horizontal line with labels $1,2,\ldots,n$. For a block $B=\{i_1,i_2,\ldots,i_k\}$  with $k\geq 2$ and
$\min(B)=i_1$, draw an arc from $i_1$ to $i_j$ for  $j = 2,\ldots,k$.
For example, the linear representation of the  linked partition
$\{1,2,4\}\{2,3\}\{3,9\}\{5,6\}\{6, 7\}\{8\}$   is illustrated in Figure \ref{LP}.

For $i<j$, we  use a pair $(i,j)$ to denote an arc from $i$ to $j$,  and we call $i$  and $j$ the left-hand endpoint and  the right-hand endpoint of $(i,j)$, respectively.
Two arcs $(i_1, j_1)$ and $(i_2, j_2)$ form a crossing if $i_1 < i_2 < j_1 < j_2$, and  form a nesting if $i_1 < i_2 < j_2 < j_1$. For the linked partition in Figure \ref{LP}, there is one crossing formed by  $(1,4)$ and $(3,9)$, while there are three nestings: $(1,4)$ and $(2,3)$, $(3,9)$ and $(5,6)$, $(3,9)$ and $(6,7)$.
A linked partition is called noncrossing (resp., nonnesting)
if there  does not exist any crossing  (resp.,  nesting) in its linear representation.
 Dykema \cite{Dyke07} showed that the number of noncrossing linked partitions of $[n+1]$ is equal to the $n$-th large Schr\"{o}der number.
  The  sequence of the large Schr\"oder numbers
   is listed as A006318 in OEIS \cite{OEIS}. Chen, Wu and Yan \cite{CWY08}
 found a bijective proof of this relation and proved that the number of 
 nonnesting linked partitions of $[n]$ also equals the number of noncrossing linked partitions 
 of $[n]$.

 Permutation tableaux are introduced by Steingr\'{i}msson and Williams \cite{SteWil07} in the study of totally positive Grassmannian cells. They are closely related to the PASEP (partially asymmetric exclusion process) model in statistical physics \cite{CortBra06, CorWil071, CorWil072, CorWil11}. Permutation tableaux are also in one-to-one correspondence with alternative tableaux introduced by Viennot \cite{Vien07}. More precisely,
  a permutation tableau is defined by a Ferrers diagram possibly with empty rows such that the cells are filled with $0$'s and $1$'s, and
\begin{itemize}
\item[(1)] each column contains at least one $1$;
\item[(2)] there does not exist a $0$ with a $1$ above (in the same column) and a $1$ to the left (in the same row).
\end{itemize}

The length of a permutation tableau is defined as the number of rows plus the number of columns.
A permutation tableau $T$ of length $n$ is labeled by the elements in $[n]$ in increasing order from the top right corner to the bottom left corner.
 We use $(i,j)$ to denote the cell in row  $i$ and column $j$.
For example, Figure \ref{PT} gives a permutation tableau of length $11$ with an empty row.

\begin{figure}[h]
\setlength{\unitlength}{0.15cm}
\begin{center}
\begin{picture}(25,32)
\put(0,0){
\begin{picture}(0,0)
\put(0,0){\line(0,1){30}}
\put(5,5){\line(0,1){25}}\put(10,10){\line(0,1){20}}
\put(15,15){\line(0,1){15}}\put(20,20){\line(0,1){10}}\put(25,20){\line(0,1){10}}
\put(0,5){\line(1,0){5}}\put(0,10){\line(1,0){10}}\put(0,15){\line(1,0){15}}
\put(0,20){\line(1,0){25}}\put(0,25){\line(1,0){25}}\put(0,30){\line(1,0){25}}
\put(25.2,26.6){\footnotesize $1$}\put(25.2,21.6){\footnotesize $2$}
\put(15.3,16.6){\footnotesize $5$}\put(10.2,11.2){\footnotesize $7$}
\put(5.3,6.4){\footnotesize $9$}\put(0.2,0.6){\footnotesize $11$}
\put(22,18){\footnotesize $3$}\put(17,18){\footnotesize $4$}
\put(12,13){\footnotesize $6$}\put(7,8){\footnotesize $8$}
\put(1.4,3){\footnotesize $10$}
\put(2,26.4){$1$}\put(2,21.4){$0$}\put(2,16.4){$0$}\put(2,11.4){$0$}\put(2,6.4){$1$}
\put(7,26.4){$0$}\put(7,21.4){$1$}\put(7,16.4){$0$}\put(7,11.4){$0$}
\put(12,26.4){$0$}\put(12,21.4){$0$}\put(12,16.4){$1$}
\put(17,26.4){$1$}\put(17,21.4){$1$}
\put(22,26.4){$0$}\put(22,21.4){$1$}
\end{picture}
}
\end{picture}
\caption{A permutation tableau  of length $11$.}\label{PT}
\end{center}
\end{figure}

 Corteel and Nadeau \cite{CorNad09} gave a bijection from  permutation tableaux of length $n$ with $k$ columns and permutations of $[n]$ with $k$ descents. Steingr\'{i}msson and Williams \cite{SteWil07} established a one-to-one correspondence between  permutation tableaux of length $n$ with $k$ rows and  permutations of $[n]$ with $k$ weak excedances.

The aim of this paper is to  demonstrate that
linked partitions play a role
as an intermediate structure between permutations and
permutation tableaux. More precisely, we
present two bijections.
The first is  between linked partitions and
permutations, and the second is between linked partitions and permutation tableaux.
 In fact, the first bijection maps a linked partition of
$[n]$ with $k$ blocks to a permutation on $[n]$ with $k-1$ descents,  and
the second bijection transforms a linked partition of $[n]$ with $k$ blocks to
a permutation tableau of length $n$ with $k$ rows.

Combining the above two bijections, we are led to a one-to-one correspondence between permutations and permutation tableaux, which implies some known  properties of permutation tableaux obtained by Corteel and Nadeau \cite{CorNad09} and Steingr\'{i}msson and Williams \cite{SteWil07}.

 When restricting the second bijection to noncrossing linked partitions, we find that
the corresponding permutation tableaux are exactly those that avoid a pattern $J_2$. Similarly, when restricting this bijection to nonnesting linked partitions,
we get  permutation tableaux that avoid a pattern $I_2$.
The definitions of the patterns $I_2$ and $J_2$ will be given in Section 4.


\section{Linked partitions and permutations}\label{SeclpP}

In this section, we   give a bijection between linked partitions of $[n]$
with $k$ blocks and permutations of $[n]$ with $k-1$ descents.

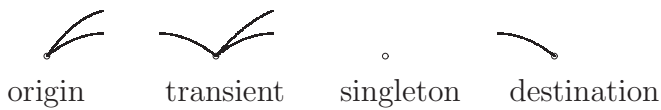
\begin{figure}[h]
\begin{center}
\setlength{\unitlength}{0.15cm}
\begin{picture}(45,10)
\put(0,-1.5){
\multiput(0,5)(15, 0){4}{\circle{0.5}}
\put(-3.5,1){origin}\put(10.5,1){transient}\put(26,1){singleton}\put(41,1){destination}
\qbezier(0,5)(2.5,7)(5,7)\qbezier(0,5)(2.5,8.5)(5,9)
\qbezier(15,5)(17.5,7)(20,7)\qbezier(15,5)(17.5,8)(20,9)\qbezier(15,5)(12.5,7)(10,7)
\qbezier(45,5)(42.5,7)(40,7)
}
\end{picture}
\caption{Four types of vertices in linked partitions.}\label{LPsort}
\end{center}
\end{figure}

To describe the construction of this bijection, we give a classification of the vertices
in the linear representation of a linked partition.
Let $\tau$ be a linked partition of $[n]$. A vertex
 $i$  in
the linear representation of $\tau$ is called an origin if it is only a left-hand endpoint, or
  a transient if it is both a left-hand point and a right-hand endpoint, or
  a singleton if
it is an isolated vertex, or a destination if it is only a right-hand endpoint.
Figure \ref{LPsort} illustrates the four types of vertices.

Let $L(n,k)$ denote the set of linked partitions of $[n]$ with $k$ blocks. Let $\pi=\pi_1 \pi_2\cdots \pi_n$ be a
 permutation on $[n]$. An integer $i$ $(1\leq i\leq n-1)$ is called a decent (resp., ascent) of $\pi$ if $\pi_i>\pi_{i+1}$ (resp., $\pi_i<\pi_{i+1}$). Let $P(n,k)$ denote the set of  permutations of $[n]$ with $k$ descents,
 which is counted by
the Eulerian number $A(n,k+1)$, see,
  for example, Stanley \cite{Stanbook}.

\begin{theo}\label{th2}
There is a bijection between $L(n,k)$ and $P(n,k-1)$.
\end{theo}

\pf We  construct a bijection $\varphi$ between $L(n,k)$ and $P(n,k-1)$ by a recursive
procedure. Let $\tau\in L(n,k)$ be a linked partition of $[n]$ with $k$ blocks. It is easily seen that the total number of  origins, transients and singletons equals $k$.
Let $\pi$ denote the permutation $\varphi(\tau)$, which
is given as follows.

If $n=1$, that is, $\tau=\{1\}$, then we set $\pi=1$.
We now assume that  $n\geq 2$. Let $\tau'$ be the linked partition of $[n-1]$ obtained by removing the vertex $n$ along with the arcs associated to $n$ from $\tau$. Assume  that $\tau'$ has $s$ blocks. Let
$i_1< i_2<\cdots< i_s$
be the minimal elements of the  blocks in $\tau'$. Let $j_1<j_2<\cdots <j_t$ be the destinations of $\tau'$. Clearly,  $t=n-1-s$.
Let  $\pi'=\pi_1\pi_2\cdots \pi_{n-1}=\varphi(\tau')$.
We  assume that the number of descents of $\pi'$  is one less than  the number of blocks in $\tau'$.

We proceed to construct a permutation $\pi$ by inserting $n$ into $\pi'$ such that the number of descents of $\pi$ is one less than the number of blocks of $\tau$.
Let $k$ be the minimal element in the block of $\tau$ containing $n$. Here are four cases.

\noindent
Case 1:
$k=i_r$, where $1\leq r\leq s-1$, that is, there is an arc $(i_r,n)$ in the linear representation of $\tau$. Set
 \begin{equation}\label{eq1}\pi=\pi_1 \pi_2 \cdots \pi_\ell  \ n \ \pi_{\ell+1} \cdots \pi_{n-1},\end{equation}
 where $\ell$ is the $r$-th descent of $\pi'$ from left to right.

 \noindent
Case 2: $k=i_s$. Set $\pi=\pi'\ n$.

\noindent
Case 3:
$k=j_r$, where $1\leq r \leq t$.  Set
\[\pi=\pi_1 \pi_2 \cdots  \pi_\ell\ n\ \pi_{\ell+1} \cdots  \pi_{n-1},\]
 where $\ell$ is the $r$-th ascent in $\pi'$ from left to right.

\noindent
Case 4: $k=n$, that is, $n$ is a singleton in $\tau$. Set
 $\pi=n\ \pi'$.

In any of the above cases, it can be shown that
$\pi$ is a permutation in  $P(n,k-1)$.
We shall only consider Case 1, since similar arguments apply to other cases.
In Case 1, it is easy to check that  $\tau'$ belongs to $L(n-1,k)$. Hence $\pi'$ belongs to $P(n-1,k-1)$.
By \eqref{eq1}, we find that $\pi$
 has the same number of descents as $\pi'$. It follows that
 $\pi$ belongs to $P(n,k-1)$.

It is straightforward to verify that
the above procedure is reversible. Hence the map $\varphi$
is a bijection for any $n\geq 1$. This completes the proof. \qed

For example, Figure \ref{eg1} illustrates the linked partitions of $\{1,2,3\}$ and the  corresponding permutations.

\begin{figure}[h]
\setlength{\unitlength}{0.8cm}
\begin{center}
\begin{picture}(18.5,3)
\qbezier(-.4,-0.5)(-0.4,1)(-0.4,2.7)\qbezier(2.6,-0.5)(2.6,1)(2.6,2.7)
\qbezier(5.8,-0.5)(5.8,1)(5.8,2.7)\qbezier(9,-0.5)(9,1)(9,2.7)
\qbezier(12.2,-0.5)(12.2,1)(12.2,2.7)\qbezier(15.4,-0.5)(15.4,1)(15.4,2.7)
\qbezier(18.6,-0.5)(18.6,1)(18.6,2.7)

\put(-.4,-.5){\line(1,0){19}}
\put(-.4,2.7){\line(1,0){19}}
\put(-.4,0.4){\line(1,0){19}}

\put(0,0){\begin{picture}(0,0)
\put(0.6,-0.2){$\footnotesize 123$}
\multiput(0,1)(1,0){3}{\circle{0.15}}
\put(-0.1,0.6){$\scriptstyle 1$}\put(0.9,0.6){$\scriptstyle 2$}\put(1.9,0.6){$\scriptstyle 3$}
\qbezier(0,1)(0.5,1.5)(1,1)\qbezier(0,1)(1,1.9)(2,1)
\put(0.1,1.8){$\footnotesize \{1,2,3\}$}

   \end{picture}
}
\put(3.2,0){\begin{picture}(0,0)
\put(0.6,-0.2){$\footnotesize 132$}
\multiput(0,1)(1,0){3}{\circle{0.15}}
\put(-0.1,0.6){$\scriptstyle 1$}\put(0.9,0.6){$\scriptstyle 2$}\put(1.9,0.6){$\scriptstyle 3$}
\qbezier(0,1)(0.5,1.6)(1,1)\qbezier(1,1)(1.5,1.6)(2,1)
\put(-0.2,1.8){$\footnotesize \{1,2\}\{2,3\}$}

   \end{picture}
}

\put(6.4,0){\begin{picture}(0,0)
\put(0.6,-0.2){$\footnotesize 231$}
\multiput(0,1)(1,0){3}{\circle{0.15}}
\put(-0.1,0.6){$\scriptstyle 1$}\put(0.9,0.6){$\scriptstyle 2$}\put(1.9,0.6){$\scriptstyle 3$}
\qbezier(0,1)(1,1.9)(2,1)
\put(0,1.8){$\footnotesize  \{1,3\}\{2\}$}

   \end{picture}
}

\put(9.6,0){\begin{picture}(0,0)
\put(0.6,-0.2){$\footnotesize 312$}
\multiput(0,1)(1,0){3}{\circle{0.15}}
\put(-0.1,0.6){$\scriptstyle 1$}\put(0.9,0.6){$\scriptstyle 2$}\put(1.9,0.6){$\scriptstyle 3$}
\qbezier(0,1)(.5,1.6)(1,1)
\put(0,1.8){$\footnotesize  \{1,2\}\{3\}$}

   \end{picture}
}
\put(12.8,0){\begin{picture}(0,0)
\put(0.6,-0.2){$\footnotesize 213$}
\multiput(0,1)(1,0){3}{\circle{0.15}}
\put(-0.1,0.6){$\scriptstyle 1$}\put(0.9,0.6){$\scriptstyle 2$}\put(1.9,0.6){$\scriptstyle 3$}
\qbezier(1,1)(1.5,1.6)(2,1)
\put(0,1.8){$\footnotesize  \{1\}\{2,3\}$}

   \end{picture}
}
\put(16,0){\begin{picture}(0,0)
\put(0.6,-0.2){$\footnotesize 321$}
\multiput(0,1)(1,0){3}{\circle{0.15}}
\put(-0.1,0.6){$\scriptstyle 1$}\put(0.9,0.6){$\scriptstyle 2$}\put(1.9,0.6){$\scriptstyle 3$}
\put(-.2,1.8){$\footnotesize  \{1\}\{2\}\{3\}$}

   \end{picture}
}
\end{picture}
\end{center}
\caption{Linked partitions of $\{1,2,3\}$ and the corresponding permutations. }\label{eg1}
\end{figure}

\section{Linked partitions and permutation tableaux}\label{section3}

The objective of this section is to
 give a bijection between
 linked partitions of $[n]$ with $k$ blocks and
 permutation tableaux of length $n$ with $k$ rows.
As consequences, we find some equidistribution
properties between  linked partitions
and permutation tableaux with certain restrictions.
For example, the number of permutation tableaux of length $n$ with $k$ rightmost restricted $0$'s is equal to the number of linked partitions of $[n]$ with $k$ transients.
We also show that our bijections can be used to
deduce some equidistribution
properties of permutation tableaux obtained by Corteel and Nadeau \cite{CorNad09} and Steingr\'{i}msson and Williams \cite{SteWil07}.

To construct the   bijection between linked partitions  and
 permutation tableaux,
we introduce the  shape of a linked partition.
Let $\tau$ be a linked partition of $[n]$.
The shape of $\tau$
is an integer partition
$\lambda$
of $n$ with empty parts allowed that is defined as follows.
Recalling the labeling of a permutation tableau, we assign labels on the boundary of  a partition $\lambda$ of $n$ in increasing order from the top right corner to the bottom left corner, see Figure \ref{PT}.
Based on  this labeling,
we see that $\lambda$ can be represented by a sequence
of labels $V$ and $H$ starting with $V$, where $V$ stands for
the vertical direction and $H$ stands for the horizontal direction.

\begin{figure}[h]
\begin{center}
\setlength{\unitlength}{0.15cm}
\begin{picture}(15,35)
\put(0,0){\line(0,1){30}}\put(5,0){\line(0,1){30}}\put(10,5){\line(0,1){25}}
\put(15,15){\line(0,1){15}}
\put(0,0){\line(1,0){5}}\put(0,5){\line(1,0){10}}\put(0,10){\line(1,0){10}}
\put(0,15){\line(1,0){15}}\put(0,20){\line(1,0){15}}\put(0,25){\line(1,0){15}}
\put(0,30){\line(1,0){15}}
\end{picture}
\caption{The shape of $\{1,2,4\}\{2,3\}\{3,9\}\{5,6\}\{6, 7\}\{8\}$.}\label{taushape}
\end{center}
\end{figure}
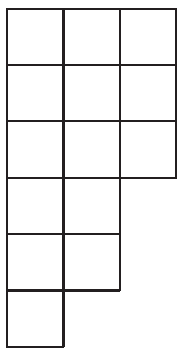

Using the  linear representation of $\tau$,
we obtain a sequence $S=s_1s_2\cdots s_n$ of
directions. If $i$ is a destination, we set $s_i=H$; otherwise, we set $s_i=V$. For example, Figure \ref{taushape}
gives the shape  of  the linked partition
in Figure \ref{LP}.

The construction of the  bijection
between linked partitions and permutation tableaux
also requires a fact proved by Corteel and Nadeau \cite{CorNad09} that a permutation tableau is determined by its topmost $1$'s and rightmost restricted
$0$'s. A $1$ in  a permutation tableau is called  a topmost $1$ if there is no $1$ above it in the same column.
A $0$ is said to be restricted if there is a $1$ above.
The rightmost such $0$ in a row is called a rightmost restricted $0$. For example, in the permutation tableau in Figure \ref{PT}, the topmost $1$'s are in the cells $(1,4)$, $(1,10)$, $(2,3)$, $(2,8)$ and $(5,6)$, while the rightmost restricted $0$'s are in the  cells $(2,10)$, $(5,8)$ and $(7,8)$.
To see that a permutation tableau is uniquely determined by its topmost $1$'s   and rightmost restricted $0$'s,
 it suffices to observe the fact
 that if  the positions of  topmost $1$'s are
 given, then all the cells above them are filled with $0$'s; if  the positions of the rightmost restricted $0$'s are given, then all the cells to the left of them are filled with $0$'s; the rest of the cells are filled with $1$'s.

Let $T(n,k)$ denote the set of permutation tableaux of length $n$ with $k$ rows. Then we have the following
correspondence and  the explicit construction is given the proof.

\begin{theo}\label{th3}
There is  a shape preserving bijection  between  $L(n,k)$ and $T(n,k)$.
\end{theo}

\pf We  construct a shape preserving bijection $\phi$ between $L(n,k)$ and $T(n,k)$.
Let $\tau$ be a linked partition in $L(n,k)$. We proceed
to construct a permutation tableau $T=\phi(\tau)$.

First, we generate the shape $\lambda$ of  $\tau$, and label the boundary of $\lambda$
by using  the elements of $[n]$ in increasing order from the top right corner to the bottom left corner.

Next, we wish to fill the cells of $\lambda$ with topmost $1$'s and rightmost restricted $0$'s.
Let $i_1$ be the minimum  origin in $\tau$, and let $j$ be a destination such that there exists a path $(i_1,i_2,\ldots,i_m, j)$ from $i_1$ to $j$ with $j$ being maximum in the linear representation of $\tau$. Then, we fill the cell $(i_1,j)$ with $1$. For $\ell=2,\ldots,m$,
we fill the cells  $(i_\ell,j)$ with $0$.

Let $\tau'$ be the linked partition of $[n]$ obtained from $\tau$
by removing the arcs in the path $(i_1,i_2,\ldots,i_m, j)$.
Repeating the  above process  for $\tau'$, we can fill a column with a topmost $1$ and some rightmost restricted $0$'s  until there  are no arcs left in the linear representation of $\tau$.

Finally, we define $T$ to be the permutation tableau such that the  $1$'s
and $0$'s
filled in $\lambda$ constitute  the topmost 1's and rightmost restricted 0's of $T$, respectively.
As mentioned before, a permutation tableau is uniquely
determined by its topmost 1's and rightmost restricted 0's.

 The reverse map of $\phi$ can be described as follows. Let $T$ be a permutation tableau of length $n$ and of shape $\lambda$. We can construct a linked partition of $[n]$. We start with the leftmost column. For the column labeled with $j$, let $(i,j)$ be the cell filled with a topmost $1$. If there does not exist any rightmost restricted $0$ in column $j$, then let $(i,j)$ be an arc in the linear representation of $\tau$. Otherwise, let $(i_2,j), (i_3,j),\ldots, (i_m,j)$ be the cells filled with the rightmost restricted $0$'s.
For  $\ell=1,2,\ldots,m$,   let  $(i_{\ell}, i_{\ell+1})$ be the arc  in the linear representation of $\tau$, where $i_1=i$ and $i_{m+1}=j$.
 Note that there is a unique topmost $1$ in each column and a unique  rightmost restricted $0$ in each row. Thus, for any vertex $j$ in $\tau$, there is at most one arc whose right-hand endpoint is  $j$. This implies that $\tau$ is a linked partition of $[n]$. Moreover, it is not hard to check that $T$ and $\tau$ have the same shape.
 This completes the proof. \qed

\begin{figure}[h]
\begin{center}
\setlength{\unitlength}{0.15cm}
\begin{picture}(75,40)
\put(0,2){\begin{picture}(0,0)
\put(2,-2.2){\footnotesize $9$}\put(5.4,1.6){\footnotesize $8$}
\put(7,2.8){\footnotesize $7$}\put(10.4,6.6){\footnotesize $6$}
\put(10.4,11.6){\footnotesize $5$}\put(12,12.8){\footnotesize $4$}
\put(15.4,16.6){\footnotesize $3$}\put(15.4,21.6){\footnotesize $2$}
\put(15.4,26.6){\footnotesize $1$}
\put(0,0){\line(0,1){30}}\put(5,0){\line(0,1){30}}\put(10,5){\line(0,1){25}}
\put(15,15){\line(0,1){15}}
\put(0,0){\line(1,0){5}}\put(0,5){\line(1,0){10}}\put(0,10){\line(1,0){10}}
\put(0,15){\line(1,0){15}}\put(0,20){\line(1,0){15}}\put(0,25){\line(1,0){15}}
\put(0,30){\line(1,0){15}}
\put(20,15){\vector(1,0){7}}
   \end{picture}
}
\put(30,2){\begin{picture}(0,0)
\put(2,-2.2){\footnotesize $9$}\put(5.4,1.6){\footnotesize $8$}
\put(7,2.8){\footnotesize $7$}\put(10.4,6.6){\footnotesize $6$}
\put(10.4,11.6){\footnotesize $5$}\put(12,12.8){\footnotesize $4$}
\put(15.4,16.6){\footnotesize $3$}\put(15.4,21.6){\footnotesize $2$}
\put(15.4,26.6){\footnotesize $1$}
\put(0,0){\line(0,1){30}}\put(5,0){\line(0,1){30}}\put(10,5){\line(0,1){25}}
\put(15,15){\line(0,1){15}}
\put(0,0){\line(1,0){5}}\put(0,5){\line(1,0){10}}\put(0,10){\line(1,0){10}}
\put(0,15){\line(1,0){15}}\put(0,20){\line(1,0){15}}\put(0,25){\line(1,0){15}}
\put(0,30){\line(1,0){15}}
\put(20,15){\vector(1,0){7}}
\put(7,6.5){$0$}
\put(7,11.5){$1$}
\put(2,16.5){$0$}
\put(2,21.5){$0$}
\put(2,26.5){$1$}
\put(12,26.5){$1$}
   \end{picture}
}
\put(60,2){\begin{picture}(0,0)
\put(2,-2.2){\footnotesize $9$}\put(5.4,1.6){\footnotesize $8$}
\put(7,2.8){\footnotesize $7$}\put(10.4,6.6){\footnotesize $6$}
\put(10.3,11.6){\footnotesize $5$}\put(12,12.8){\footnotesize $4$}
\put(15.4,16.6){\footnotesize $3$}\put(15.4,21.6){\footnotesize $2$}
\put(15.4,26.6){\footnotesize $1$}
\put(0,0){\line(0,1){30}}\put(5,0){\line(0,1){30}}\put(10,5){\line(0,1){25}}
\put(15,15){\line(0,1){15}}
\put(0,0){\line(1,0){5}}\put(0,5){\line(1,0){10}}\put(0,10){\line(1,0){10}}
\put(0,15){\line(1,0){15}}\put(0,20){\line(1,0){15}}\put(0,25){\line(1,0){15}}
\put(0,30){\line(1,0){15}}
\put(7,6.5){$0$}
\put(7,11.5){$1$}
\put(2,16.5){$0$}
\put(2,21.5){$0$}
\put(2,26.5){$1$}
\put(12,26.5){$1$}
\put(2,1.5){$1$}\put(2,6.5){$0$}
\put(2,11.5){$1$}\put(7,16.5){$0$}\put(7,21.5){$0$}\put(7,26.5){$0$}
\put(12,16.5){$1$}\put(12,21.5){$1$}
   \end{picture}
}
\end{picture}
\caption{From a linked partition to a permutation tableau.}\label{phitau}
\end{center}
\end{figure}
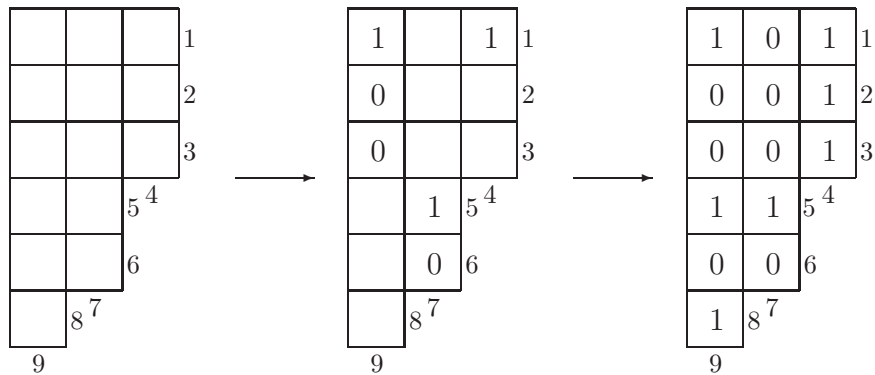

 For example,
 the permutation tableau corresponding to the linked partition in Figure \ref{LP} is given in Figure \ref{phitau}.

\begin{figure}[h]
\begin{center}
\setlength{\unitlength}{0.7cm}
\begin{picture}(10,2.7)
\multiput(0,.7)(1,0){11}{\circle{0.1}}
\qbezier(0,.7)(.5,1.2)(1,.7)\qbezier(1,.7)(1.5,1.2)(2,.7)
\qbezier(0,.7)(1.5,2.1)(3,.7)\qbezier(1,.7)(2.5,1.9)(4,.7)
\qbezier(1,.7)(5,3.7)(9,.7)\qbezier(4,.7)(4.5,1.2)(5,.7)
\qbezier(4,.7)(5,1.75)(6,.7)\qbezier(6,.7)(6.5,1.2)(7,.7)
\put(-0.1,0){$1$}\put(.9,0){$2$}\put(1.9,0){$3$}\put(2.9,0){$4$}\put(3.9,0){$5$}
\put(4.8,0){$6$}\put(5.9,0){$7$}\put(6.9,0){$8$}\put(7.8,0){$9$}\put(8.75,0){$10$}
\put(9.75,0){$11$}
\end{picture}
\caption{The  linked partition corresponding to the permutation tableau in Figure \ref{PT}.}\label{PTlp}
\end{center}
\end{figure}
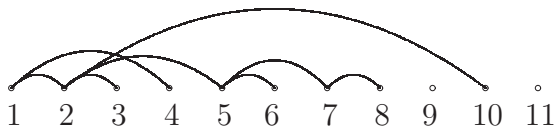

For the permutation tableau given in Figure \ref{PT}, the corresponding linked partition is given in Figure \ref{PTlp}.

It is easy to see that
the bijection
$\phi$ has the following properties.

\begin{coro} For $n\geq 1$, let $\tau$ be a
linked partition of $[n]$. Then the number of arcs in the linear representation of $\tau$ is equal to the total number of topmost $1$'s and rightmost restricted $0$'s in $T=\phi(\tau)$.
\end{coro}

\begin{coro} Assume that $n\geq 1$.
\begin{itemize}
\item[(1)] For $0\leq k \leq n-2$, the number of linked partitions of $[n]$ with $k$ transients equals the number of permutation tableaux of length $n$ with $k$ rightmost restricted $0$'s;
\item[(2)]
For $0\leq k \leq n$, the number of linked partitions of $[n]$ with $k$ singletons equals the number of permutation tableaux of length $n$ with $k$ rows that do not contain any topmost $1$ or restricted $0$.
\end{itemize}
\end{coro}

To conclude this section, we remark that our bijections $\varphi$ and $\phi$ can be used to
deduce
some results on permutation tableaux obtained by  Corteel and Nadeau \cite{CorNad09}
and
Steingr\'{i}msson and Williams \cite{SteWil07}.
 Corteel and Nadeau \cite{CorNad09} showed that the number of permutation tableaux
 of length $n$ with $k$ columns is equal to the number of permutations of $[n]$ with $k$ descents.
This fact follows from bijections $\varphi$ and $\phi$. Noting  that the number of permutations of $n$ with $k-1$ descents is equal to the number of permutations
 of $[n]$ with $n-k$ descents, we see that the number of permutation tableaux of length $n$ with $k$ rows is equal to
the number of permutations with $n-k$ descents, which is equivalent to the result of
 Corteel and Nadeau \cite{CorNad09}.

On the other hand, the number of permutations of $[n]$
with $k-1$ descents is equal to the number of permutations of $[n]$
with $k$ weak excedances, see, for example, Stanley \cite[Chapter 1]{Stanbook}.
This leads to the fact  that the number of permutation tableaux
of length $n$ with $k$ rows is equal to the
number of permutations of $[n]$ with $k$ weak  excedances,
as proved by Steingr\'{i}msson and Williams \cite{SteWil07}.

\section{Pattern avoiding permutation tableaux}

In this section, we discuss   restrictions of the  bijection $\phi$ in Section \ref{section3} to noncrossing linked partitions
and nonnesting linked partitions, and characterize
the corresponding permutation tableaux by pattern avoidance.

\begin{figure}[h]
\setlength{\unitlength}{0.13cm}
\begin{center}
\begin{picture}(20,27)

\put(3,7){\thicklines
\begin{picture}(0,0)
\hspace*{-3cm}
\put(0,0){\line(0,1){17}}\put(5,0){\line(0,1){17}}\put(14,0){\line(0,1){17}}
\put(19,0){\line(0,1){17}}
\put(0,0){\line(1,0){19}}\put(0,5){\line(1,0){19}}\put(0,12){\line(1,0){19}}
\put(0,17){\line(1,0){19}}
\put(2.4,14.2){\circle*{1.1}}\put(16.4,2.2){\circle*{1.1}}
\put(8,-7.3){$I_2$}
\put(2.3,7){$\vdots$}\put(16.2,7){$\vdots$}\put(7.5,8.1){$\cdots$}
\put(7.5,2.1){$\cdots$}\put(7.5,14.3){$\cdots$}
\put(0.8,-2.8){$j_2$}\put(16,-2.8){$j_1$}\put(19.4,1){$i_2$}\put(19.4,14){$i_1$}
\hspace*{5cm}

\put(0,0){\line(0,1){17}}\put(5,0){\line(0,1){17}}\put(14,0){\line(0,1){17}}
\put(19,0){\line(0,1){17}}
\put(0,0){\line(1,0){19}}\put(0,5){\line(1,0){19}}\put(0,12){\line(1,0){19}}
\put(0,17){\line(1,0){19}}
\put(16.4,14.2){\circle*{1.1}}
\put(2.4,2.2){\circle*{1.1}}
\put(8,-7.3){$J_2$}
\put(2.3,7){$\vdots$}\put(16.2,7){$\vdots$}\put(7.5,8.1){$\cdots$}
\put(7.5,2.1){$\cdots$}\put(7.5,14.3){$\cdots$}
\put(0.8,-2.8){$j_2$}\put(16,-2.8){$j_1$}\put(19.4,1){$i_2$}
\put(19.4,14){$i_1$}
\end{picture}
}

\end{picture}
\caption{Patterns $I_2$ and $J_2$.}\label{pattern}
\end{center}
\end{figure}
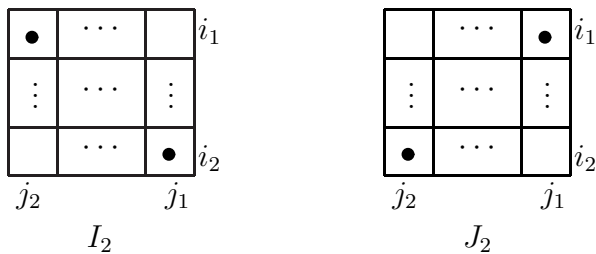

We introduce two patterns $I_2$ and $J_2$ as given in Figure \ref{pattern},
where a dot means a topmost $1$ or a rightmost restricted $0$.
 We use $I_2$ and $J_2$ to
 denote these two patterns because similar
notation  has been used in the context of fillings of Ferrers diagrams, see de Mier \cite{DeM07}.

\begin{figure}[h]
\setlength{\unitlength}{0.12cm}
\begin{center}
\begin{picture}(50,38)
\put(3,3){
\begin{picture}(0,0)\thicklines
\put(0,0){\line(0,1){30}}\put(5,0){\line(0,1){30}}\put(10,0){\line(0,1){30}}
\put(15,0){\line(0,1){30}}\put(20,0){\line(0,1){30}}\put(25,10){\line(0,1){20}}
\put(30,10){\line(0,1){20}}\put(35,20){\line(0,1){10}}
\put(0,0){\line(1,0){20}}\put(0,5){\line(1,0){20}}\put(0,10){\line(1,0){30}}
\put(0,15){\line(1,0){30}}\put(0,20){\line(1,0){35}}\put(0,25){\line(1,0){35}}
\put(0,30){\line(1,0){35}}
\put(35.2,26.6){\footnotesize $1$}\put(35.2,21.6){\footnotesize $2$}
\put(30.1,16.5){\footnotesize $4$}\put(30.2,11.3){\footnotesize $5$}
\put(20.2,6.6){\footnotesize $8$}\put(20.2,0.8){\footnotesize $9$}
\put(32,17.8){\footnotesize $3$}\put(27,7.8){\footnotesize $6$}
\put(22,7.8){\footnotesize $7$}\put(16.2,-2.2){\footnotesize $10$}
\put(11.2,-2.2){\footnotesize $11$}\put(6.2,-2.2){\footnotesize $12$}
\put(1.2,-2.2){\footnotesize $13$}
\put(2.4,27.2){\circle*{1.1}}\put(2.4,7.2){\circle*{1.1}}
\put(2.4,22.2){\circle*{1.1}}
\put(7.4,7.2){\circle*{1.1}}
\put(12.4,7.2){\circle*{1.1}}
\put(17.4,7.2){\circle*{1.1}}
\put(22.4,22.2){\circle*{1.1}}
\put(27.4,22.2){\circle*{1.1}}
\put(32.4,27.2){\circle*{1.1}}
\end{picture}
}
\end{picture}
\caption{A permutation tableau avoiding  $I_2$.}\label{POA}
\end{center}
\end{figure}
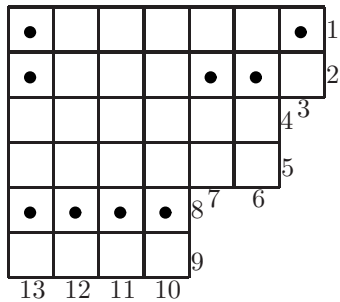

More precisely, let $T$ be a permutation tableau of length $n$, and let $T'$ be the
permutation tableau obtained from $T$ by replacing the topmost $1$'s and rightmost restricted $0$'s by dots and removing all other $1$'s and $0$'s.
We say that $T$ avoids the pattern $I_2$ if $T'$ does not contain  four cells $(i_1,j_1)$, $(i_1,j_2)$, $(i_2,j_1)$ and $(i_2,j_2)$, where $i_1<i_2<j_1<j_2$, such that  the cells $(i_1,j_2)$ and $(i_2,j_1)$ are filled with dots, while the cell $(i_2,j_2)$ is empty.
Similarly, we say that $T$ avoids the pattern $J_2$ if $T'$ does not contain four cells $(i_1,j_1)$, $(i_1,j_2)$, $(i_2,j_1)$ and $(i_2,j_2)$, where $i_1<i_2<j_1<j_2$,  such that
 the cells $(i_1,j_1)$ and $(i_2,j_2)$ are filled with dots, while the cell $(i_2,j_1)$ is empty.
For example, Figure \ref{POA} illustrates a permutation tableau avoiding the pattern $I_2$.
Note that this permutation tableau does not avoid the pattern $J_2$,  because the cells $(1,3)$,
$(1,6)$, $(2,3)$ and $(2,6)$  form   pattern $J_2$.

By the construction of the bijection $\phi$ in Section \ref{section3},
we obtain   characterizations of
permutation tableaux corresponding to noncrossing
linked partitions and nonnesting linked partitions.
To be more specific, we have the following correspondences.

\begin{theo}
There is a bijection between noncrossing linked partitions of $[n]$ and
$J_2$-avoiding permutation tableaux of length $n$, and there is a bijection
between nonnesting linked partitions of $[n]$ and $I_2$-avoiding permutation tableaux of length $n$.

\end{theo}

Since  the number of noncrossing linked partitions of $[n]$ equals the number of 
nonnesting linked partitions of $[n]$, see Chen, Wu and Yan \cite{CWY08},
 one sees that the number of 
 $J_2$-avoiding permutation tableaux of length $n$ equals the number of
 $I_2$-avoiding permutation tableaux of length $n$.

\vspace{15pt}

\noindent  {\bf Acknowledgments.} This work was supported by the 973 Project, the PCSIRT Project of the Ministry of Education, the National Science Foundation of China, Beijing Natural Science Foundation, and Beijing Commission of Education.

\end{document}